\documentclass[12pt]{article}
	\usepackage{pp}
	\usepackage{ppmacro}
	\usepackage[colorlinks=true,%
		linkcolor=blue,%
		citecolor=blue]{hyperref}
	\usepackage[dvips]{graphicx}

	\newcommand{\chn}{\check{n}}
	\newcommand{\cC}{\mathcal{C}}
	\newcommand{\cH}{\mathcal{H}}
	\newcommand{\cL}{\mathcal{L}}
	\newcommand{\cK}{\mathcal{K}}
	
	\newcommand{\cT}{\mathcal{T}}
	\newcommand{\fsl}{\mathfrak{sl}}
	\newcommand{\fg}{\mathfrak{g}}
	
	\newcommand{\fz}{\mathfrak{z}}
	\newcommand{\Zcs}{Z_{CS}}
	\newcommand{\vmu}{\vec{\mu}}
	\newcommand{\vA}{\vec{A}}
	\newcommand{\vB}{\vec{B}}
	\newcommand{\vx}{\vec{x}}
	\newcommand{\vy}{\vec{y}}
	
	\newcommand{\cY}{\mathcal{Y}}

	\newcommand{\bZ}{\mathbb{Z}}
	
	\newcommand{\qnum}[1]{ \big( q^{-\frac{#1}{2} }- q^{\frac{#1}{2} } \big) }
	
	\newcommand{\mathcenter}[1]{ \vcenter{ \hbox{#1} } }

\begin{document}
	\title{ New structures of knot invariants }
			\author{Kefeng Liu and Pan Peng}
			\date{}
	\maketitle

	\begin{abstract}
		Based on the proof of Labastida-Mari{\~n}o-Ooguri-Vafa conjecture \cite{lmov}, we
	derive an infinite product formula for Chern-Simons partition functions, the generating
	function of quantum $\fsl_N$ invariants. Some symmetry properties of the infinite product
	will also be discussed.
	\end{abstract}


\section{Introduction}
	Chern-Simons theory has been conjectured to be equivalent to a topological string theory by
	$1/N$ expansion in physics. This duality conjecture builds a fundamental connection in
	mathematics. On the one hand, Chern-Simons theory leads to the
	construction of knot invariants; on the other hand, topological string theory gives rise to
	Gromov-Witten theory in geometry.

	Therefore, the Chern-Simons/topological string duality conjecture identifies the generating
	function of Gromov-Witten invariants as Chern-Simons knot invariants \cite{OV}.  Based on these
	thoughts, the existence of a sequence of integer invariants is conjectured \cite{OV, LMV, MV, LM}
	in a similar spirit of Gopakumar-Vafa setting \cite{GV}, which provides an essential evidence of
	the duality between Chern-Simons theory and topological string theory.
	This integrality conjecture, called the LMOV conjecture, was proved in \cite{lmov}.

	One important corollary of the LMOV conjecture is to express Chern-Simons partition function as an
	infinite product derived in this article. The motivation of studying such an infinite product formula is
	based on a guess on the modularity property of topological string partition function.

	Chern-Simons theory can be approached in mathematics with the help of quantum group
	theory. Let $\cL$ be a link of $L$ components. Quantum $\fsl_N$ invariant of $\cL$,
	$W_{(A^1,\ldots,A^L)}(\cL;q,t)$, is defined to be a trace function on the $U_q(\fsl_N)$ modules
	constructed from the planar diagram of $\cL$ and irreducible $U_q(\fsl_N)$ modules
	$V_{A^1},\ldots, V_{A^L}$ associated to the components of $\cL$. These irreducible
	$U_q(\fsl_N)$ modules are labeled by Young diagrams $A^1,\ldots, A^L$.
	However, we will define these invariants more explicitly by the HOMFLY skein model
	described in    section \ref{section 3}.

	The Chern-Simons partition function of $\cL$ is the following generating function of quantum
	$\fsl_N$ invariants:
	\begin{align*}
		\Zcs(\cL; q,t) = 1+ \sum_{A^1,\ldots, A^L} W_{(A^1,\ldots, A^L)}(\cL; q,t)
			\prod_{\alpha=1}^L s_{A^\alpha}(x^\alpha),
	\end{align*}
	where the summation is taken over all the partitions with some of them possibly being empty while
	not all empty, and $s_{A^\alpha}(x^\alpha)$ is the Schur function of a set of variables
	$x^\alpha=\{ x^\alpha_i\}_{i\geq 1}$. Based on our proof of LMOV conjecture,
	we obtain the following infinite product formula (we only show the infinite product
	formula for a knot. Its generalization to links is similar. Please refer to section \ref{section 4}) :
	\begin{align*}
		\Zcs = \prod_\mu 
			\prod_{Q\in \bZ/2}\; \prod_{m=1}^\infty \; \prod_{k=-\infty}^\infty\;					
					\big \langle 1- q^{k+m}t^Q  x^\mu
				\big \rangle^{-m\check{n}_{\mu;\,k,Q}}\,.
	\end{align*}
	Here, $\langle \cdot \rangle$ is the symmetric product defined by
	\eqref{definition of symmetric product}, and $\check{n}_{\mu;k,Q}$ are invariants related to
	the integer invariants in the LMOV conjecture. The symmetry property of $\chn_{\mu;k,Q}$ is
	discussed in section \ref{section 5}. For more details of the infinite product formula for
	the Chern-Simons partition function of a link $\cL$,
	please refer to \eqref{infinite product formula for a link}.

	The paper is organized as follows. In section \ref{section 2}, we define the quantum $\fsl_N$
	invariants by the HOMFLY skein model. Chern-Simons theory and LMOV conjecture are described
	in section \ref{section 3}. In section \ref{section 4}, we will derive Chern-Simons partition
	function as infinite product. In section \ref{section 5}, we will discuss  the symmetry of 
	$q\rightarrow q^{-1}$ and the rank-level duality.


\section{Quantum $\fsl_N$ invariants}\label{section 2}
\subsection{Preliminary}
	We start by reviewing some preliminary on the representations of symmetric groups and
	symmetric functions.

	A partition  of $n$ is a tuple of
	positive integers $\mu=(\mu_1, \mu_2, \ldots , \mu_k )$ such that
	$|\mu|\triangleq\sum_{i=1}^k \mu_i =n$
	and $\mu_1 \geq \mu_2\geq \cdots \geq \mu_k>0$, where $|\mu|$ is called the degree
	of $\mu$ and $k$ is called the length of $\mu$, denoted by $\ell(\mu)$.
	A partition can be represented by a Young diagram,
	for example, partition $(5,4,2,1)$ can be identified as the following Young diagram:
	\begin{align*}
		\partition{5 4 2 1} \,.
	\end{align*}
	Denote by $\cY$ the set of all Young diagrams.
	Let $\chi_A$ be the character of irreducible representation of
	symmetric group, labelled by partition $A$ . Given a partition $\mu$,
	define $m_j = \textrm{card} (\mu_k=j; k\geq 1)$.
	The order of the conjugate class of type $\mu$ is given by:
	\begin{align*}
		\fz_\mu = \prod_{j\geq1} j^{m_j} m_j!.
	\end{align*}
	The orthogonality of the character formula gives
	\begin{align*}
		\sum_{\mu} \frac{ \chi_A(C_\mu) \chi_B(C_\mu) }{ \fz_\mu } = \delta_{A,B}
		= \bigg \{
			\begin{array}{ll}
				1, & \text{if } A=B; \\
				0, & \text{otherwise}.
			\end{array}
	\end{align*}
	The theory of symmetric functions has a close relationship with the representations of
	symmetric group. The
	symmetric power functions of a given set of variables $x=\{x_j\}_{j\geq 1}$
	are defined as the direct limit of the Newton polynomials:
	\begin{align*}
		p_n(x) = \sum_{j\geq1} x_j^n,\qquad p_\mu(x) = \prod_{i\geq 1} p_{\mu_i}(x),
	\end{align*}
	and we have the following formula which determines the Schur function
	\begin{align*}
		s_A(x) = \sum_{\mu} \frac{ \chi_A(C_\mu) }{\fz_\mu} p_\mu(x)\,.
	\end{align*}
	Given $x= \{x_i\}_{i\geq 1}$, $y=\{y_j\}_{j\geq 1}$, define
	\begin{align}\label{product of two sets of invariables}
		x\ast y = \{x_i\cdot y_j\}_{i\geq 1, j\geq 1}.
	\end{align}	
	We also define $x^d = \{ x_i^d\}_{i\geq 1}$.
	The $d$-th Adam operation of a Schur function is given by $s_A(x^d)$.

\subsection{The HOMFLY skien models}
	The quantum group invariants can be defined over any semi-simple Lie algebra $\fg$. Here,
	we are particularly interested in the $SU(N)$ Chern-Simons gauge theory, hence we will study
	the quantum $\fsl_N$ invariants, which can be identified as the colored HOMFLY
	polynomials.

	We will start by introducing the Framed HOMFLY polynomial, an invariant of framed oriented links.
	Define the skein $R_n^n(q,t)$ by linear combinations of oriented $n$-tangles modulo the following
	relations:
	\begin{align}\label{skein relation}
	\begin{array}{rcl}
		\mathcenter{ \includegraphics[height=1cm]{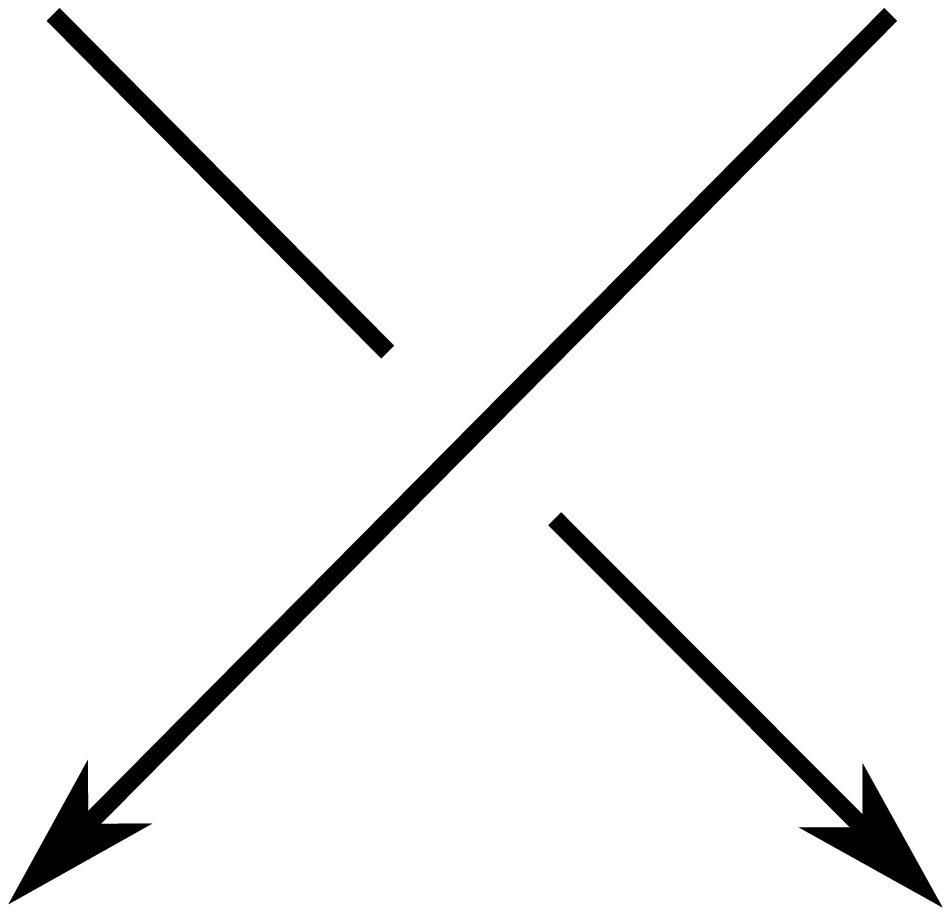}
					}
		- \ \  \mathcenter{\includegraphics[height=1cm]{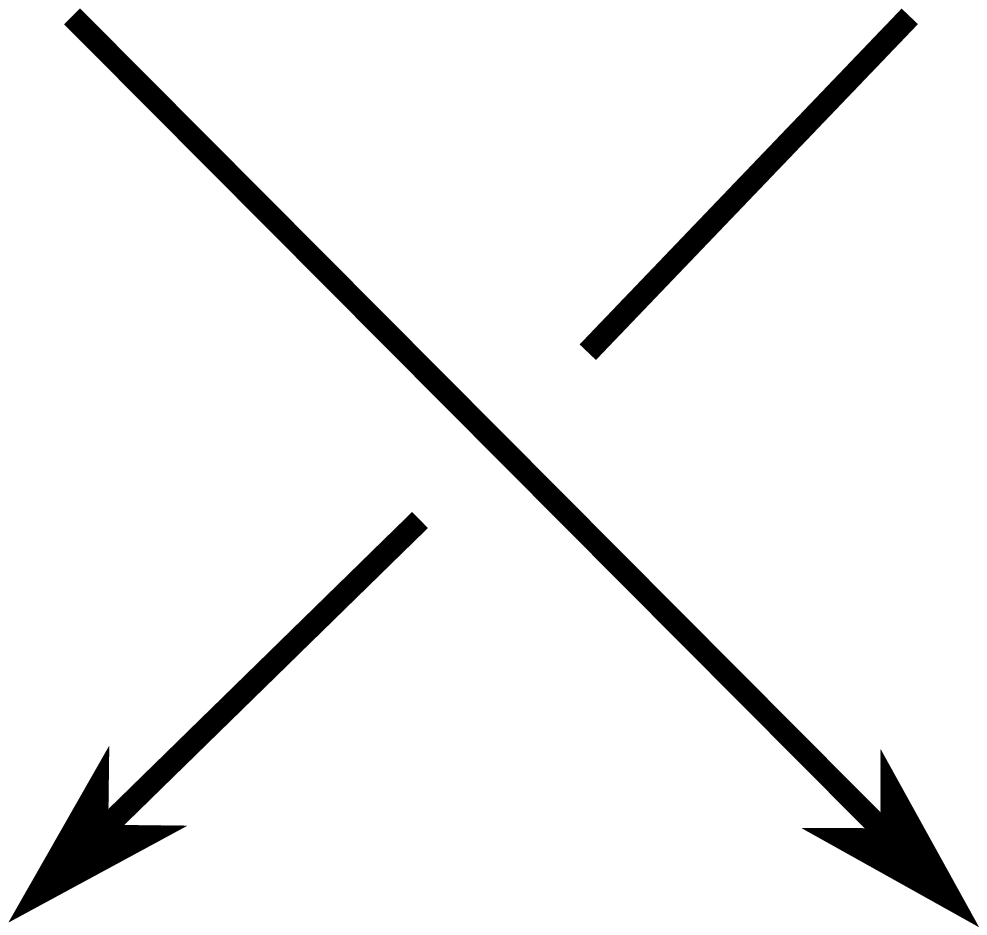} } \
		&=&\ \qnum{1}\ \  \mathcenter{\includegraphics[height=1cm]{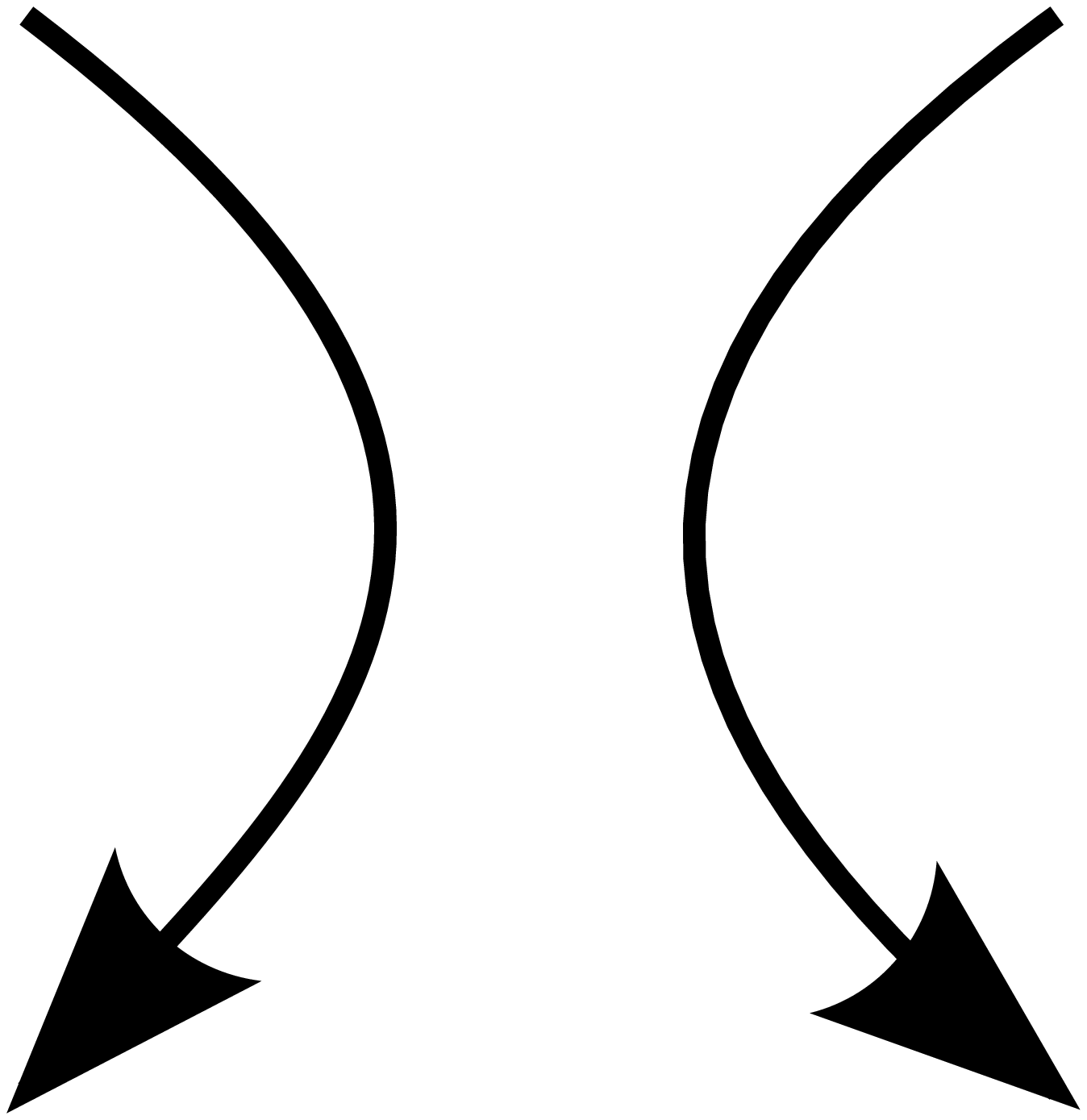}  }\,,
	\\
		&&
	\\
		\textrm{and} \qquad
		\mathcenter{ \includegraphics[height=1.2cm] {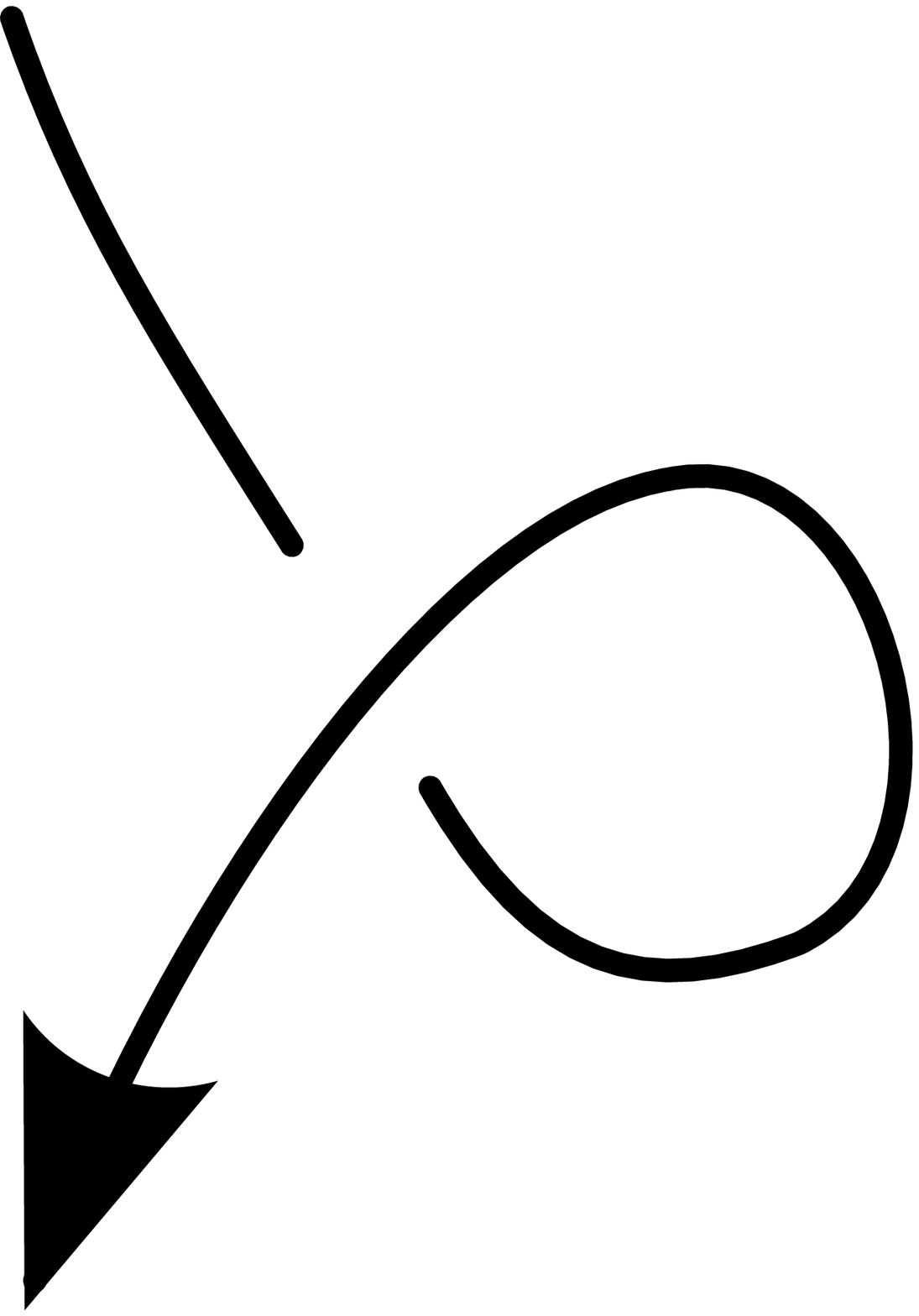} }  \
		&=& \ t^{-\frac12} \ \mathcenter{\includegraphics[height=1.2cm]{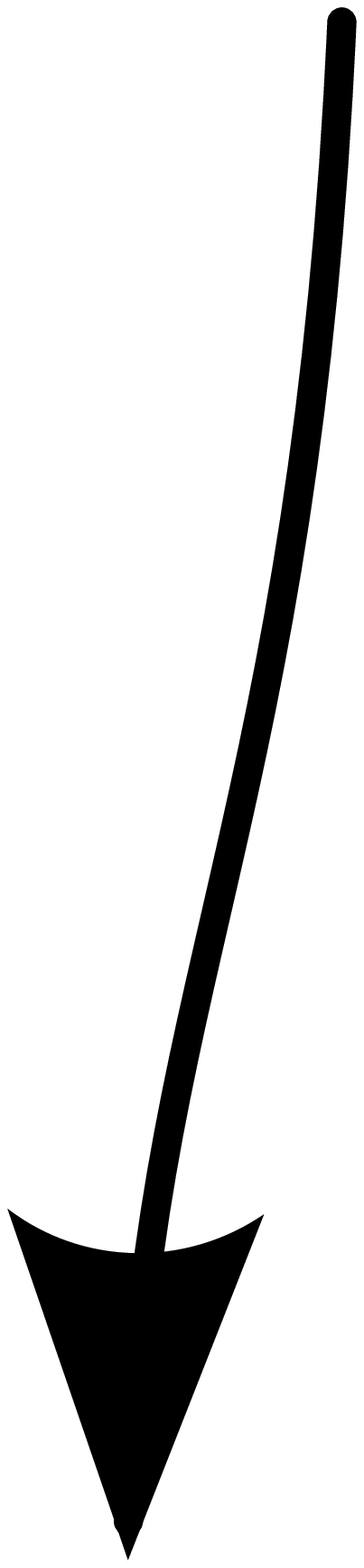}     } \,.
	\end{array}
	\end{align}
	The HOMFLY skein of an planar annulus, with some designated input and output boundary points,
	is defined as linear combinations of oriented tangles in the annulus, modulo Reidemeister moves II
	and III and the above two relations \eqref{skein relation}. The coefficient ring is
	$\mathbb{Z} \big[ q^{\pm\frac 12}, t^{\pm\frac12} \big]$ with finitely many products of
	$\qnum{k}$ in the denominators. The skein of the annulus is denoted by $\cC$. This is a commutative
	algebra with a product given by placing one annulus outside another.

	From the HOMFLY skein of annulus, we can obtain the framed HOMFLY polynomial of links, denoted
	by $\cH(\cL)$. Here we normalize $\cH$ as:
	\begin{align*}
		\cH(\mathrm{unknot}) = \frac{ t^{-\frac12}-t^{\frac12} }{ q^{-\frac12}-q^{\frac12} }\,.
	\end{align*}
	These invariants can be recursively computed through the HOMFLY skein.

\subsection{The quantum group invariants}
	The colored HOMFLY polynomials are defined through \emph{satellite knot}. A satellite of $\cK$ is
	determined by choosing a diagram $Q$ in the annulus. Draw $Q$ on the annular neighborhood of
	$\cK$ determined by the framing to give a satellite knot $\cK\star Q$.
	See the following figure for a satellite of a framed trefoil knot with $Q$:
	\begin{align*}
		\cK = \mathcenter{ \includegraphics[height=1in]{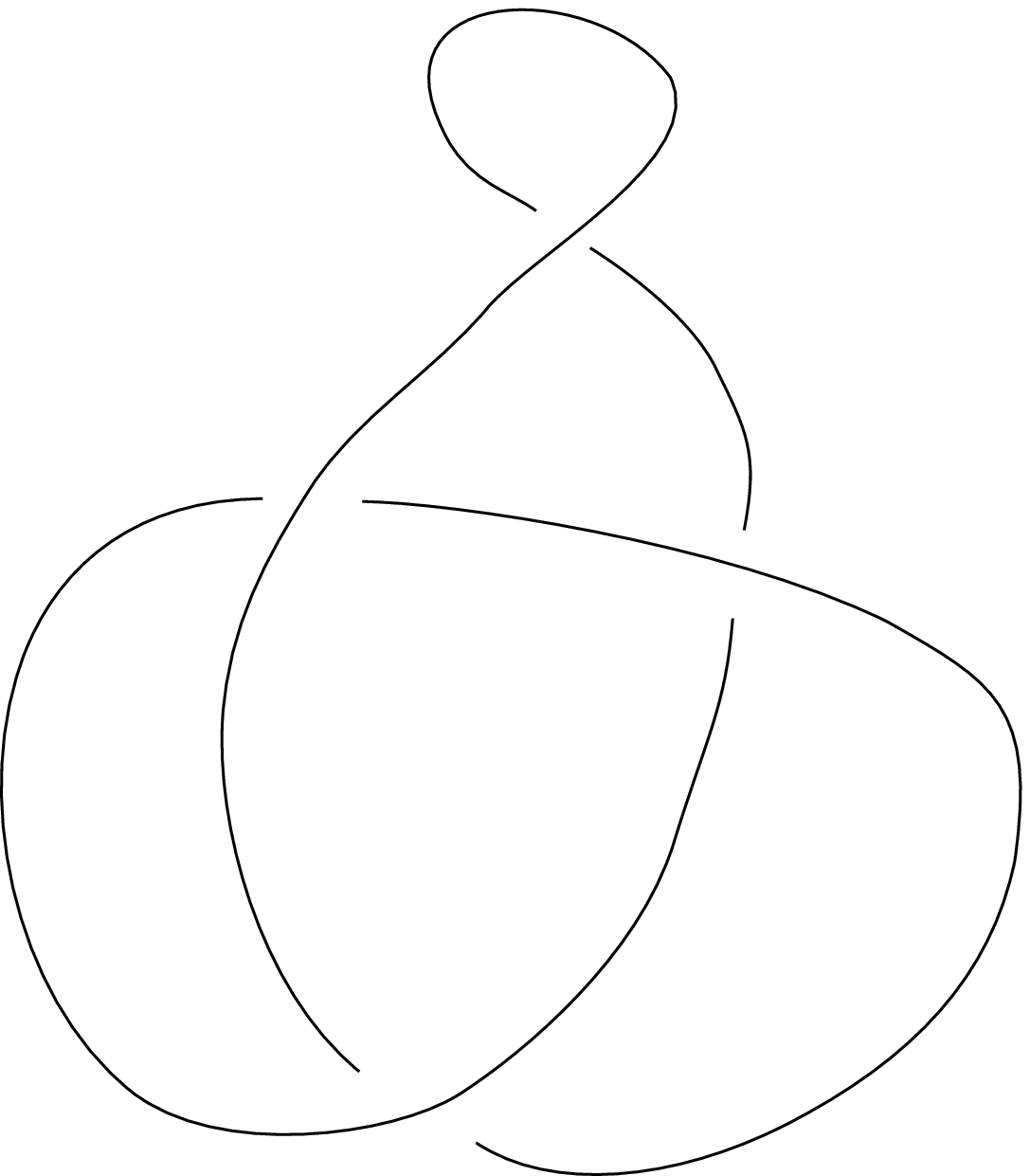}
									},
		&&
		Q = \mathcenter{ \includegraphics[height=0.8in]{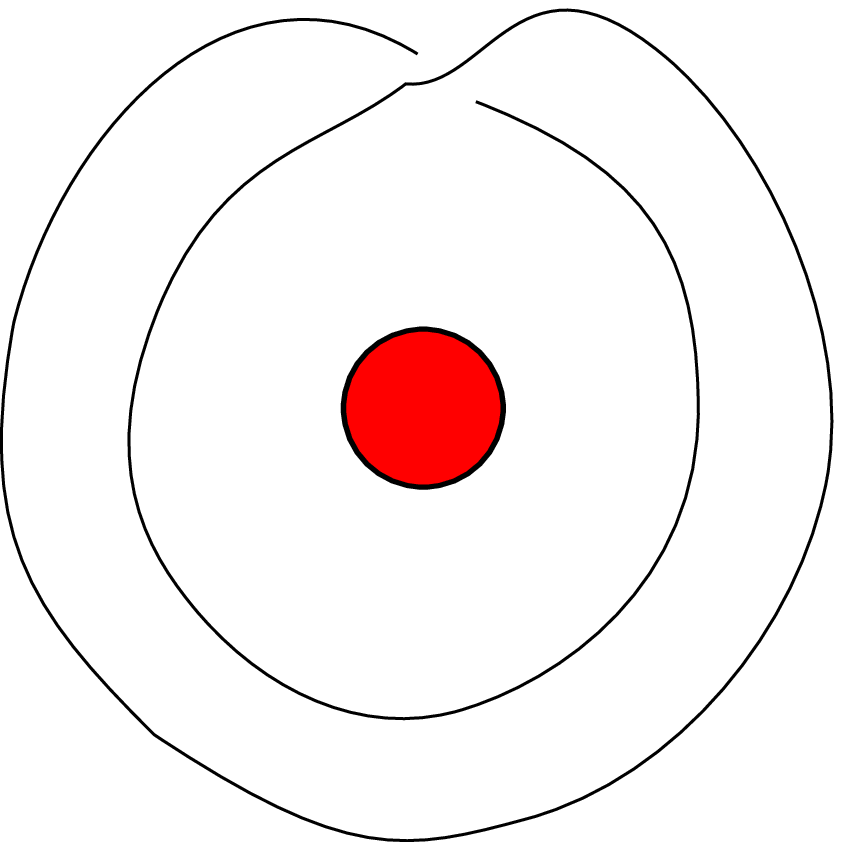}
									},
		&&
		\cK\star Q = \mathcenter{ \includegraphics[height=1in]{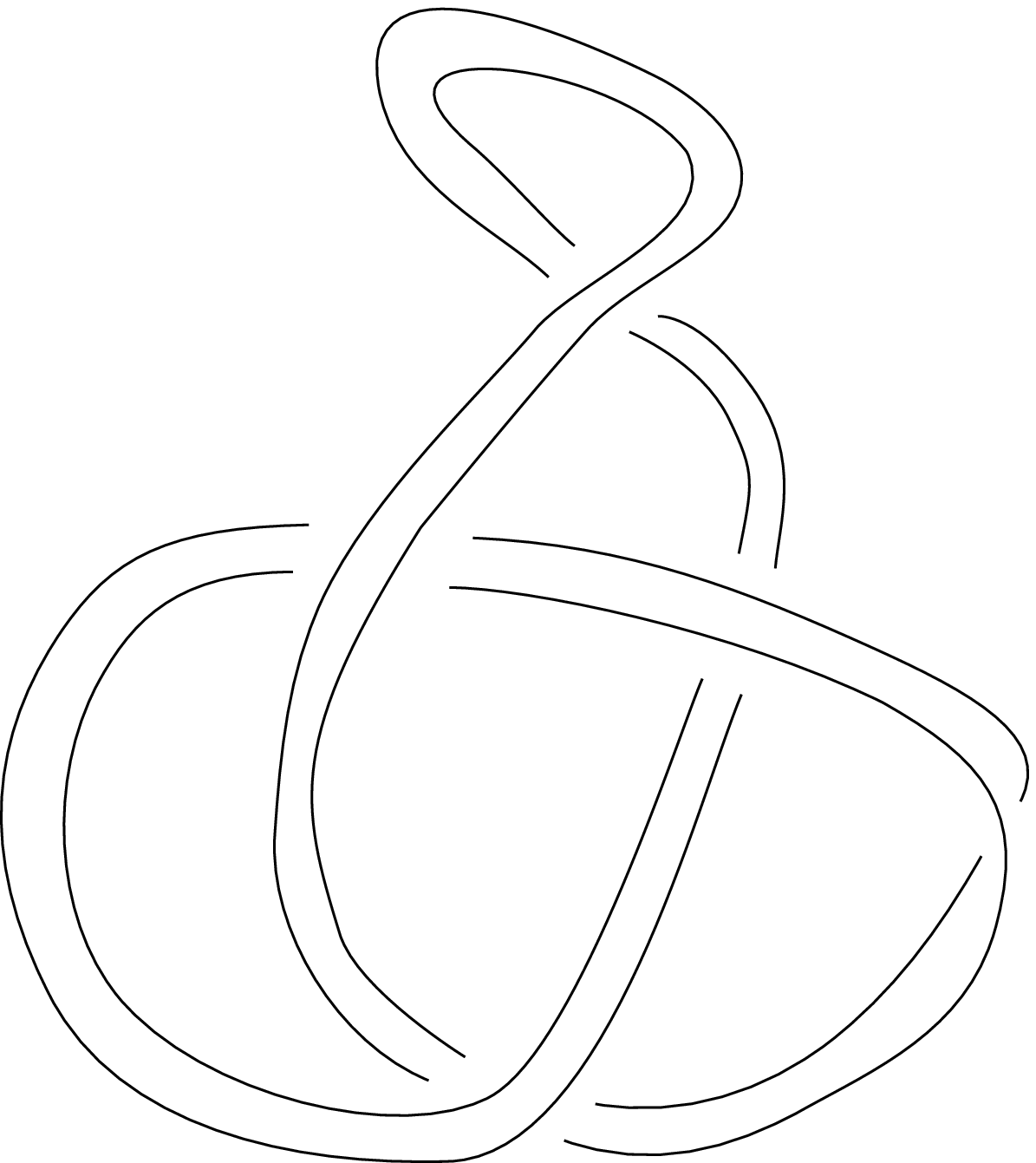}
											}.
	\end{align*}

	We will refer to this construction as \emph{decorating $\cK$ with the pattern $Q$}. The HOMFLY
	polynomial $\cH(\cK\star Q)$ of the satellite depends on $Q$ only as an element of the skein
	$\cC$ of the annulus. $\cC$ can be regarded as the parameter space for these invariants of
	$\cK$. We will call it the \emph{HOMFLY satellite invariants of $\cK$}.

	There is a known set of idempotent elements, $E_\lambda$, one for each partition $\lambda$ of
	$n$. They were originally described in \cite{Gyoja, Morton-Aiston}.
	Take the closure of $E_\lambda$, we have $Q_\lambda\in\cC$.
	$\{Q_\lambda\}_{\lambda\in\cY}$ form a basis of $\cC$.

	The quantum $\fsl_N$ invariant for the irreducible module $V_{A^1},\ldots,V_{A^L}$,
	labeled by the corresponding partitions $A^1,\ldots, A^L$,  can be
	identified as the HOMFLY invariants for the link decorated by $Q_{A^1},\ldots,Q_{A^L}$.
	Write $\vA=(A^1,\ldots,A^L)$. The quantum $\fsl_N$ invariants of the link is given by
	\begin{align}\label{definition of quantum sl_N invariants}
		W_{\vA}(\cL; q,t) = \cH (\cL\star \otimes_{\alpha=1}^L Q_{A^\alpha} )\,,
	\end{align}
	and we can define the following invariants:
	\begin{align}\label{definition of Z_mu}
		Z_{\vec{\mu}}(\cL;q,t) = \sum_{\vec{A}=(A^1,\ldots,A^L) } \bigg( \prod_{\alpha=1}^L
			\chi_{A^\alpha}(C_{\mu^\alpha} ) \bigg)
			W_{\vec{A} }(\cL;q,t)\,.
	\end{align}

\subsection{Notations}
	Here we make the following convention for the notations in this article.
	\begin{itemize}
		\item We will consistently denote by $\cL$ a link and by $L$ the number of
			components in $\cL$.
		\item The irreducible $U_q(\fsl_N)$ module associated to $\cL$ will be labelled by
			their highest   weights, thus by Young diagrams. We usually denote it by a vector form
			$\vec{A}=(A^1,\ldots,A^L)$.
		\item Let $\vx=(x_1,\ldots,x_L)$ is $L$ sets of variables,
			each of which is associated to  a component of $\cL$ and
			$\vmu = (\mu^1,\ldots,\mu^L)\in\cY^L$ be a tuple of $L$ partitions. Define:
			\begin{align*}
				&[n] = q^{-\frac n2}-q^{\frac n2}, &    &[\vmu] = \prod_{\alpha=1}^L [\mu^\alpha], &
					& \fz_{\vmu} = \prod_{\alpha=1}^L \fz_{\mu^\alpha},
				\\
				&\chi_{\vA}(C_{\vmu}) = \prod_{\alpha=1}^L \chi_{A^\alpha}(C_{\mu^\alpha}), &
					& s_{\vA}(\vx) = \prod_{\alpha=1}^L s_{A^\alpha}(x_\alpha), &
					& p_{\vmu}(\vx) = \prod_{\alpha=1}^L p_{\mu^\alpha}(x_\alpha).
			\end{align*}
	\end{itemize}


\section{Chern-Simons theory}\label{section 3}
\subsection{Chern-Simons partition function}
	The \emph{Chern-Simons partition function} of the link $\cL$ is the following generating series of
	quantum group invariants weighted by Schur functions:
	\begin{align*}
		\Zcs(\cL;q,t) = 1+ \sum_{\vec{A}} W_{\vec{A}} (\cL; q,t) s_{\vec{A}}(\vx) \,,
	\end{align*}
	or equivalently by the reformulated invariants $Z_{\vec{\mu}}$:
	\begin{align}\label{partition in terms of Z_mu}
		\Zcs(\cL;q,t) = 1+ \sum_{\vec{\mu} } \frac{ Z_{\vmu} (\cL;q,t) }{\fz_{\vmu} } p_{\vmu}(\vx) \,.
	\end{align}
	\emph{Free energy} is defined by
	\begin{align}\label{free energy}
		F(\cL;q,t) = \log \Zcs(\cL;q,t)
		= \sum_{\vmu} \frac{ F_{\vmu}(\cL;q,t) }{\fz_{\vmu}} p_{\vmu}(\vx) \,.
	\end{align}
	Here, we expand the free energy and get the definition of $F_{\vmu}(\cL;q,t)$ according to
	the above formula.

\subsection{Transformation function}
	Define the following \emph{transformation function} for two partitions, $A$ and $B$, of $n$:
	\begin{align}\label{definition of transformation function}
		\cT_{AB}(x) = \sum_{\mu} \frac{ \chi_A(C_\mu) \chi_B(C_\mu) } {\fz_\mu} p_\mu(x)\,.
	\end{align}
	By a simple algebra, its inverse is given by
	\begin{align}\label{definition of the inverse of transformation function}
		\cT^{-1}_{AB}(x) = \sum_{\mu} \frac{ \chi_A(C_\mu) \chi_B(C_\mu) } {\fz_\mu}
			\frac{1}{ p_\mu(x) }  \,.
	\end{align}
	One immediately sees the following:
	\begin{align}\label{transformation formula for schur function}
		\sum_{A} \cT_{AB}(x) s_A(y)
		&=\sum_A \sum_{\mu} \frac{ \chi_A(C_\mu) \chi_B(C_\mu) }{\fz_\mu} p_\mu(x) s_A(y) \notag\\
		&= \sum_{\mu} \frac{ \chi_B(C_\mu) }{\fz_\mu} p_\mu(y) \sum_A \chi_A(C_\mu) s_A(x)\notag\\
		&=  \sum_{\mu} \frac{ \chi_B(C_\mu) }{\fz_\mu} p_\mu(y) p_\mu(x) \notag \\
		&= s_B ( x\ast y )\,.
	\end{align}
	Here in the last step, $x\ast y$ is defined in \eqref{product of two sets of invariables}.
	Similarly, we have
	\begin{align}\label{inverse transformation formula for schur function}
		\sum_A \cT^{-1}_{AB} (x) s_A (x\ast y) = s_B(y)
	\end{align}
	Define variables $q^\varrho = \{ q^{j-\frac12} \}_{j\geq 1}$ and
	$q^{d\varrho}= \{ q^{d(j-\frac12)} \}_{j\geq 1}$. We have
	\begin{align}\label{q-number as symmetric power function}
		\frac{1}{[n]} = p_n(q^\varrho)\,.
	\end{align}

\subsection{LMOV conjecture}
	Based on the large $N$ duality between Chern-Simons gauge theory and topological string
	theory, Labastida, Mari{\~n}o, Ooguri, and Vafa made a conjecture on the existence of a series of
	integer invariants in Chern-Simons theory \cite{OV,LMV,LM}.

	Free energy has the following expansion
	\begin{align*}
		F = \sum_{\vA} \sum_{d=1}^\infty \frac 1d f_{\vA}(q^d, t^d) s_{\vA}\big( (\vx)^d \big)\,,
	\end{align*}
	where
	\begin{align*}
		s_{\vA}\big( (\vx)^d \big) = \prod_{\alpha=1}^L s_{A^\alpha}
			\Big( \{(x^\alpha_j)^d\}_{j\geq 1} \Big)
	\end{align*}

	The Labastida-Mari{\~n}o-Ooguri-Vafa conjecture proved in \cite{lmov} can be
	stated as the following theorem:
	\begin{theorem}[Liu-Peng]
		Notations are as above. Define
		\begin{align}\label{definition of P_B}
			P_{\vB}(q,t) = \sum_{\vA} f_{\vA}(q,t) \prod_{\alpha=1}^L
				\cT^{-1}_{A^\alpha B^\alpha}(q^\varrho )\,.
		\end{align}
		Then
		\begin{align}\label{liu-peng theorem}
			P_{\vB}(q,t) \in [1]^{-2} \cdot\mathbb{Z}\big[ [1]^2, t^{\pm\frac12} \big].
		\end{align}
	\end{theorem}


\section{Infinite product formula}\label{section 4}
	To derive an infinite product formula, we will state the result for a knot at first, since the
	notations in the computation for a knot are relatively simpler.
\subsection{The case of a knot}\label{subsection: the case of a knot}
	Set $y= x\ast q^\varrho$, then $p_n(x\ast q^\varrho) = p_n(x)\cdot p_n(q^\varrho)$. We have
	\begin{align*}
		p_\mu(y) = p_\mu(x)p_\mu(q^\varrho).
	\end{align*}
	Consider free energy weighted by the Schur function of $y$, LMOV conjecture implies the following
	reformulation of free energy:
	\begin{align*}
		F(q,t; y) &= \sum_{d= 1}^\infty \sum_{A} \frac{1}{d} f_A(q^d,t^d) s_A(y^d)\\
		  &= \sum_{d= 1}^\infty \sum_A \frac{1}{d}
			\sum_B \cT^{-1}_{AB}(q^{d\varrho}) P_B(q^d,t^d) s_A(y^d)\,,
	\end{align*}
	where
	\begin{align}
		[1]^2\cdot P_B (q,t) &= \sum_{Q\in \bZ/2}\; \sum_{g= 0}^\infty N_{B;\,g,Q}\qnum{1}^{2g} t^Q
			\label{eqn: formular of P_B}
	\end{align}
	is a polynomial of $[1]^2$ and $t^{\pm\frac12}$ with integer coefficients, $N_{B;g,Q}$.
	By \eqref{inverse transformation formula for schur function}, we have
	\begin{align*}
		F(q,t;y) = \sum_{d=1}^\infty \sum_B \frac{1}{d} P_B(q^d,t^d) s_B(x^d)\,.
	\end{align*}
	There exist integers $n_{B;\,g,Q}$ such that
	\begin{align}\label{definition of n_B,k,Q}
		\sum_{g=0}^\infty N_{B;\,g,Q} \qnum{1}^{2g}
			= \sum_{g=0}^\infty n_{B;\,g,Q} \, \sum_{k=0}^g q^{g-2k}
		\,.
	\end{align}
	This is due to the equivalence of two integral bases of $U_q(\fsl_2)$ modules: one is obtained by
	all the irreducible  modules of $U_q(\fsl_2)$, $V_n$, for each $n\geq 0$;
	the other one is obtained by $V_1^{\otimes n}$ for $n\geq 0$ since
	$V_1\otimes V_n = V_{n+1}\oplus V_{n-1}$.

	By \eqref{liu-peng theorem} $N_{B;g,Q}$ vanish for sufficiently large $g$ and $|Q|$,
	thus $n_{B;\,g,Q}$ vanish for sufficiently large $g$ and $|Q|$.
	We have
	\begin{align*}
		P_B(q,t) &= \sum_{Q\in\bZ/2}\; \sum_{g=0}^\infty N_{B;\,g,Q} \qnum{1}^{2g-2} t^Q \\
		&= \sum_Q \frac{t^Q}{\qnum{1}^2} \sum_{g=0}^\infty n_{B;\,g,Q} \sum_{k=0}^g q^{g-2k}\\
		&= \sum_Q t^Q \sum_{m=1}^\infty m q^{m}
				\sum_{g=0}^\infty n_{B;\,g,Q} \sum_{k=0}^g q^{g-2k} \\
		&= \sum_{Q\in\bZ/2}\,\, \sum_{m=1}^\infty\,  \sum_{g=0}^\infty \sum_{k=0}^g
				m n_{B;\,g,Q}\, q^{g-2k+m} t^Q\,,
	\end{align*}
	which leads to
	\begin{align}
		F &= \sum_B \sum_{d=1}^\infty
				\sum_{Q\in\bZ/2}\,\, \sum_{m=1}^\infty\,  \sum_{g=0}^\infty \sum_{k=0}^{g} \frac{1}{d}
				m\, n_{B;\,g,Q} t^{dQ}  q^{d(g-2k+m)} s_B(x^d) \notag
	 \\
		&= \sum_{Q,m,g}\sum_{k=0}^g  \sum_{B,\mu}
				m\, n_{B;\,g,Q}   \frac{\chi_B(C_\mu)}{\fz_\mu}
				\sum_{d=1}^\infty \frac{1}{d} q^{d(g-2k+m)} p_\mu(x^d) t^{dQ} \,.
				\label{eqn: F in terms of q series}
	\end{align}
	Now focus on the following computation:
	\begin{align*}
		\sum_{d\geq1} \frac{1}{d} \omega^d p_\mu(x^d)
		&= \sum_{d\geq1} \frac{1}{d} \omega^d \prod_{j=1}^{\ell(\mu)} \sum_{i\geq1} x_k^{d\mu_j}\\
		&= \sum_{d\geq1} \frac{1}{d} \omega^d \sum_{i_1,\ldots,i_\ell} (x_{i_1}^{\mu_1} \cdots
			x_{i_\ell}^{\mu_\ell} )^d\\
		&= \sum_{i_1,\ldots, i_\ell} \sum_{d\geq1} \frac{1}{d} (\omega x_{i_1}^{\mu_1} \cdots
			x_{i_\ell}^{\mu_\ell})^d \\
		&= -\sum_{i_1,\ldots, i_\ell} \log (1-\omega x_{i_1}^{\mu_1} \cdots
			x_{i_\ell}^{\mu_\ell}).
	\end{align*}
	Let $\omega = t^Q q^{g-2k+m}$ and apply the above computation in
	\eqref{eqn: F in terms of q series},
	\begin{align*}
		F = \sum_{Q,m,k} \sum_{B,\mu} (-m\, n_{B;\,g,Q}) \frac{\chi_B(C_\mu)}{\fz_\mu}
			\sum_{i_1,\ldots,i_{\ell(\mu)} }
				\log \Big(1-q^{g-2k+m}t^Q x_{i_1}^{\mu_1}\cdots
								x_{i_{\ell(\mu)} }^{\mu_{\ell(\mu)} }
					\Big).
	\end{align*}
	Define the symmetric product as shown in the following formula:
	\begin{align}\label{definition of symmetric product}
			\big\langle 1- \omega x^\mu \big\rangle
		= \prod_{ x_{i_1},\ldots, x_{i_{\ell(\mu)} } }
			\Big( 1- \omega x_{i_1}^{\mu_1} \cdots x_{i_{\ell(\mu)}}^{\mu_{\ell(\mu)}} \Big) ,
	\end{align}
	and
	\begin{align}\label{definition of check{n}}
			\check{n}_{\mu;\,g,Q} = \sum_B \frac{\chi_B(C_\mu)}{\fz_\mu} n_{B;\,g,Q} .
	\end{align}
	Therefore, for $\Zcs=\exp F$, we obtain the following
	infinite product formula:
	\begin{align}\label{infinite product formula}
		\Zcs(\cK;q,t;y)
			  =   \prod_{\mu\in\cY}
				\prod_{Q\in \bZ/2} \,\,
				\prod_{m=1}^\infty \;
				\prod_{g=o}^\infty\;
					\prod_{k=0}^g\;
				\big \langle 1- q^{g-2k+m}t^Q   x^\mu
				\big \rangle^{-m\, \check{n}_{\mu;\,g,Q}}
	\end{align}
	\begin{remark}\label{remark of finiteness of product for k and Q}
			In the above infinite product formula, since for a given $\mu$,
		$\chn_{\mu;g,Q}$ vanish for sufficiently large
		$g$ and $|Q|$ due to the vanishing property of $n_{B;\,g,Q}$,
		the products involved with $Q$ and $g$ are
		in fact finite products for a fixed partition, $\mu$.
	\end{remark}

\subsection{The case of a link}
	After going over subsection
	\ref{subsection: the case of a knot}, we find that the computation can be carried over to the
	case of a link.
	Given a link $\cL$ of $L$ components, let $\vy=(y_1,\ldots, y_L)$ and $\vx=(x_1,\ldots,x_L)$
	satisfying $y_i= q^{\varrho}\ast x_i$, for $i=1,\cdots,L$.
	We define $n_{\vB;g,Q}$ and
	$\chn_{\vmu;g,Q}$ as the following:
	\begin{align*}
		\sum_{g=0}^\infty N_{\vB;g,Q}
	= \sum_{g=0}^\infty n_{\vB;g,Q} \sum_{k=0}^g q^{g-2k},
	\qquad
		\chn_{\vmu;g,Q}
	= \sum_{\vB} \frac{ \chi_{\vA}(C_{\vmu} ) }{ \fz_{\vmu} } n_{\vB;g,Q} \,.
	\end{align*}
	Again, the vanishing result of $N_{\vB;g,Q}$ implies that $n_{\vB;g,Q}$ and
	$\chn_{\vmu;g,Q}$ vanish for sufficiently large $g$ and $|Q|$.

	Let $\vmu=(\mu^1,\ldots,\mu^L)$ and $\vx=(x_1,\ldots, x_L)$. Denote by $\ell_i$ the length of
	$\mu^i$. Generalize the symmetric product in \eqref{definition of symmetric product} to $\vmu$
	and $\vx$ as:
	\begin{align*}
		\big\langle 1-\omega\, x_1^{\mu^1}\cdots x_L^{\mu^L}  \big\rangle
			=   \prod_{\alpha=1}^L\
				\prod_{i_{\alpha, 1},\ldots, i_{\alpha, \ell_\alpha} }
				\Big( 1- \omega
					\prod_{\alpha=1}^L
						\big(  (x_\alpha)_{ i_{\alpha,1} }^{\mu^\alpha_1} \cdots
								 (x_\alpha)_{ i_{\alpha, \ell_\alpha} }^{\mu^\alpha_{\ell_\alpha}}
						\big)
				\Big)\,.
	\end{align*}
	Follow the similar computation, we have the infinite product formula for the Chern-Simons
	partition function of $\cL$:
	\begin{align}\label{infinite product formula for a link}
		\Zcs(\cL;q,t; \vy)
			=  \prod_{\vmu\in\cY^L}\,
				\prod_{Q\in \bZ/2} \,
				\prod_{m=1}^\infty 
				\prod_{g=0}^\infty
					\prod_{k=0}^g
				\big \langle 1- q^{g-2k+m}\, t^Q \, x_1^{\mu^1}\cdots x_L^{\mu^L}
				\big \rangle^{-m\, \check{n}_{\vmu;\,g,Q}}\,.
	\end{align}
	Similar as Remark \ref{remark of finiteness of product for k and Q}, the products related to $Q$ and
	$g$ are finite for a fixed $\vmu$.

\subsection{The case of the unknot}
	The Chern-Simons partition function of the unknot is given by
	\begin{align}\label{Chern-Simons partition function of unknot}
		\Zcs(\textrm{unknot};\, q,t)
			= 1+ \sum_{A\in\cY} \dim_q V_A \cdot s_A(x).
	\end{align}
	Here, $\dim_q V_A$ is the quantum dimension of the irreducible $U_q(\fsl_N)$ module $V_A$.
	The formula of quantum dimension is given by (the formula of the quantum dimension can
	be found in many literatures, for example, \cite{lmov} ):
	\begin{align*}
		\dim_q V_A
		= \sum_\mu \frac{ \chi_A(C_\mu) }{ \fz_\mu } \prod_{j=1}^{\ell(\mu)}
			\frac{ t^{-\frac{ \mu_j} {2} } - t^{\frac{ \mu_j}{2} }  }
				{  q^{-\frac{ \mu_j} {2} } - q^{ \frac{ \mu_j}{2} } }\,.
	\end{align*}
	A similar computation as shown in the previous section leads to the following infinite product formula:
	\begin{align} \label{infinite product formula for unknot}
		\Zcs(\text{unknot};\, q,t; y)
		&= \prod_{m=1}^\infty \prod_i
			\frac{     \big( 1- q^m t^{\frac12} x_i
					\big)^{ m }
				  }%
				 {    \big( 1- q^m t^{-\frac12} x_i
					\big)^{ m }
				 }\,.%
	\end{align}
	By a comparison with \eqref{infinite product formula}, we have
	\begin{align*}
		\chn_{\mu;\, g, Q} = \delta_{\mu, \partition{1}} \cdot \delta_{g,0} \cdot \mathrm{sign}(-Q)\,.
	\end{align*}


\section{Symmetry property}\label{section 5}
	We will discuss some symmetry properties of the infinite product formula given in section
	\ref{section 4}. Here, we will express it for a knot. The case of links is exactly
	similar, otherwise we will make a remark of their difference.

\subsection{The symmetry of $q\rightarrow q^{-1}$ }
	In the derivation of the infinite product formula,
	we assume $|q|<1$ for taylor expansion of $\frac{1}{[1]^2}$. In
	the case of $|q|>1$, the taylor expansion is given by
	\begin{align*}
		\frac{1}{[1]^2} = \sum_{m=1}^\infty mq^{-m}\,.
	\end{align*}
	Therefore, the infinite product formula will be:
	\begin{align*}
		\Zcs(\cK;q,t;y)
			&=      \prod_{\mu\in\cY}
				\prod_{Q\in \bZ/2} \,\,
					\prod_{m=1}^\infty \;
					\prod_{g=o}^\infty\;
				\prod_{k=0}^g\;
					\big \langle 1- q^{g-2k}q^{-m}t^Q   x^\mu
				\big \rangle^{-m\, \check{n}_{\mu;\,g,Q}}
	\\
		&=  \prod_{\mu\in\cY}
				\prod_{Q\in \bZ/2} \,\,
					\prod_{m=1}^\infty \;
					\prod_{g=o}^\infty\;
				\prod_{k=0}^g\;
					\big \langle 1- q^{-g+2(g-k)} q^{-m}
						t^Q     x^\mu
				\big \rangle^{-m\, \check{n}_{\mu;\,g,Q}}
	\\
		&=  \prod_{\mu\in\cY}
				\prod_{Q\in \bZ/2} \,\,
					\prod_{m=1}^\infty \;
					\prod_{g=o}^\infty\;
				\prod_{k=0}^g\;
					\big \langle 1-  \big( q^{-1}  \big)^{g-2k+m}
						t^Q     x^\mu
				\big \rangle^{-m\, \check{n}_{\mu;\,g,Q}} \,.
	\end{align*}
	The substitution of $k$ to $g-k$ is used in the last equation.
	This is the symmetry of $q\rightarrow q^{-1}$ for the infinite product formula.

\subsection{Rank-level duality}
	Rank-level duality is essentially a symmetry of quantum group invariants relating a labeling
	color to its transpose, which can be expressed as the following identity:
	\begin{align}\label{rank-level duality formula}
		W_{ A^t}( s^{-1}, -v ) = W_{A} (s, v)\,,
	\end{align}
	where $s=q^{\frac12}$, $v=t^{\frac12}$.
	From \cite{lmov}, we can actually obtain the following stronger version:
	\begin{align}
		W_{A^t}( s^{-1}, v) &= (-1)^{|A|} W_A (s, v), 
			\label{strong rank-level duality, I}
	\\
		W_A(s, -v) &= (-1)^{|A|} W_A(s, v).
			\label{strong rank-level duality, II}
	\end{align}
	
	\begin{remark}
		The above two equations of the strong version rank-level duality can be interpreted by
		the $1/N$ expansion for the Chern-Simons partition function and the monodromy of 
		$t \rightarrow e^{2\pi i}t$.
	\end{remark}
	
	We now give a proof of this strong version of the rank-level duality.	
	\eqref{strong rank-level duality, I} is equivalent to proposition $5.2$ in \cite{lmov}:
	\begin{align*}
			W_{A^t}(q,t) = (-1)^{|A|} W_A(q^{-1}, t)\,.
	\end{align*}
	
	To prove \eqref{strong rank-level duality, II}, we will investigate the HOMFLY polynomial invariants
	of knots first. By the HOMFLY skein relation:
	\begin{align}\label{homfly skein relation}
		t^{-\frac12} \mathcenter{ \includegraphics[height=1cm]{L+.eps}
					}
		- \ \  t^{\frac12} \mathcenter{\includegraphics[height=1cm]{L-.eps} } \
		=\ \qnum{1}\ \  \mathcenter{\includegraphics[height=1cm]{L0.eps}  }\,,		
	\end{align}
	in the HOMFLY polynomial, the degrees of $t$ are either all integers or 
	half integers. This can be done through a resolution of the crossings. The unknot 
	obviously satisfies this property. Suppose that for any planar diagram representing a link 
	with number of crossings less or equal to $n$, the degrees of $t$ in the HOMFLY polynomial of the link
	are either integers or half integers, we immediately see that switching positive crossing to negative 
	crossing will keep this property by \eqref{homfly skein relation}. 
	This will reduce the verification to the case of disjointed unknots. 
	To generalize this result to the quantum $\fsl_N$
	invariants, we will apply (5.31) in \cite{lmov} to directly obtain that $Z_\mu(q,t)$ defined in 
	\eqref{partition in terms of Z_mu} satisfies that 
	\begin{align*}
		t^{-|\mu|/2} Z_\mu \in \mathbb{Q} \big[ \qnum{1}^{\pm1}, t^{\pm 1} \big]\,,
	\end{align*}
	by the observation that in $v^d-v^{-d}$, changing the sign of $v$ will result in a multiplication of 
	$(-1)^d$. This completes the proof of \eqref{strong rank-level duality, II}.

	\eqref{strong rank-level duality, I} leads to the symmetry of the 
	following integer invariants by expansion:
	\begin{align}\label{symmetry of N_A;g,Q}
			N_{A^t;\,g,Q} = (-1)^{|A|} N_{A;\,g, -Q }.
	\end{align}
	Applying it to \eqref{definition of check{n}} and combining the fact that
	\begin{align*}
		\chi_{A^t}(C_\mu) = (-1)^{|\mu|-\ell(\mu)} \chi_A(C_\mu)\,,
	\end{align*}
	we have the following symmetry about $\mu$ and $Q$:
	\begin{align}\label{symmetry of check{n}}
		\chn_{\mu;\, g,-Q} = (-1)^{\ell(\mu)} \chn_{\mu;\, g, Q} \,.
	\end{align}

	This is (4.44) in \cite{MV}, the rank-level duality of the $SU(N)_k$ and $SU(k)_N$ 
	Chern-Simons gauge theories.

\subsection{Concluding remarks}
	The infinite product formula derived in section \ref{section 4} is 
	interestingly related to a conjecture
	on the modularity property of the Chern-Simon partition function, 
	which is hardly seen from the knot theory point of view.  
	We hope this formula will shed a new light on the study of knot invariants. 
	A complete understanding	in this aspect will lead to a deep insight of 
	the Chern-Simons/topological string duality.




\noindent \textsc{Center of Mathematical Sciences, \\
Zhejiang University, Box 310027, \\
Hangzhou, P.R.China.}

\noindent \textsc{Department of mathematics , \\
University of California at Los Angeles, Box 951555\\ 
Los Angeles, CA, 90095-1555.} \\
Email: liu@math.ucla.edu.

\vspace{3pt}

\noindent \textsc{Department of Mathematics , \\
University of Arizona,\\
617 N. Santa Rita Ave. ,\\
Tucson, AZ, 85721.} \\
Email: ppeng@math.arizona.edu.

\end{document}